\documentclass[11pt]{article}

\usepackage{amsmath}

\def\b{\ensuremath{\binom}}
\def\s{\ensuremath{\stackrel}}
\def\la{\ensuremath{\lambda}}

\begin{document}
\begin{center}
{\Large
The Eigenvectors of the Right-Justified Pascal Triangle \\ 
}
\vspace{10mm}
DAVID CALLAN  \\
Department of Statistics  \\
University of Wisconsin-Madison  \\
1210 W. Dayton St   \\
Madison, WI \ 53706-1693  \\
{\bf callan@stat.wisc.edu}  \\
\vspace{5mm}
\today
\end{center}

\vspace{10mm}
\vspace{5mm}

Let $R = \left( \binom{i-1}{n-j} \right)_{1 \leq i,j \leq n}
$ denote the $n \times n$ matrix formed by right justifying the  
first $n$ rows of Pascal's triangle. Let $a$ denote the golden ratio 
$(1+\sqrt{5})/2$. We will show that the eigenvalues of $R$ are 
$\la_{j} = (-1)^{n+j}a^{2j-n-1},\ 1 \leq j \le n,$ (as conjectured in 
\cite{peele}), with corresponding eigenvectors 
$\mathbf{u}_{j} = (u_{ij})_{1 \le i \le n}$ where $u_{ij} = \sum_{k=1}^{j}
 (-1)^{i-k}\binom{i-1}{k-1}\binom{n-i}{j-k}a^{2k-i-1}$. 
Since the eigenvalues are distinct, the eigenvectors are linearly 
independent and so form an invertible matrix that diagonalizes $R$. 
Scaling the eigenvectors to 
$\mathbf{v}_{j} = (-1)^{j}a^{n-j}/(1+a^{2})^{(n-1)/2}\mathbf{u}_{j}$ 
yields a diagonalizing matrix $V$ 
($V^{-1}RV = \mathrm{diag}(\la_{j})_{j=1}^{n}$) with a remarkable 
property: $V^{-1} = V$. This makes it easy to write down explicit 
three-summation formulas for the entries of powers of $R$.

The proofs below that $R \mathbf{u}_{p} = \la_{p}\mathbf{u}_{p}$ for $1 \le p \le 
n$ and $V^{2} = I_{n}$ are bracing exercises in manipulating binomial 
coefficient sums. During these manipulations, summations will be 
extended over all integers when convenient; recall that a binomial 
coefficient with a negative lower parameter is zero. We must also be 
careful to avoid the symmetry trap \cite[p.\,156]{gkp}: the symmetry law 
$\binom{n}{k} = \b{n}{n-k}$ is valid only when the upper parameter $n$ 
is nonnegative.

The proof that  $R \mathbf{u}_{p} = \la_{p}\mathbf{u}_{p}$ uses
the (minimal polynomial) equation for the golden ratio: $a^{2} = a+1$, 
equivalently, $1-a = -a^{-1}$, and relies on the following two 
binomial coefficient identities.
\[
 \hspace{40mm}\b{N-J}{K} =
 \sum_{r}^{}(-1)^{r}\b{N-r}{K-r}\b{J}{r}\hspace{25mm} 
 \mathrm{integer\ }K\quad  (*)
\]
\[
 \hspace{50mm}\b{I}{J}\b{J}{K} = \b{I}{K}\b{I-K}{J-K}\hspace{25mm} 
 \mathrm{integers\ }J,K\quad  (**)
\]
The first follows from the Vandermonde convolution (below) using 
``upper negation'' to write $\binom{N-r}{K-r}$ as 
$(-1)^{K-r}\binom{K-N-1}{K-r}$, and the second is the ``trinomial 
revision'' identity \cite[p.\,174]{gkp}.
 Here goes. The \emph{i}th entry of  
$R\mathbf{u}_{p}$ is
 \begin{eqnarray*}
    & & \sum_{j=1}^{n}R_{ij}u_{jp} \\
      & = &
     \sum_{j=1}^{n}\b{i-1}{n-j}\sum_{k=1}^{p}(-1)^{j-k}
     \b{j-1}{k-1}\b{n-j}{p-k}a^{2k-j-1}  \\
    & \s{1}{=} &
    \sum_{j=1}^{n} 
    \sum_{k=1}^{p}(-1)^{n+j+p+k}\b{i-1}{j-1}\b{n-j}{p-k}\b{j-1}{k-1}
      a^{2(p+1-k)-(n+1-j)-1} \\
   &  \s{2}{=} &
   \sum_{j}^{}\sum_{k}^{}
   (-1)^{n+j+p+k}\b{i-1}{k-1}\b{i-k}{j-k}\sum_{r}^{}(-1)^{r}
   \b{n-k-r}{p-k-r}\b{j-k}{r}a^{2p-2k-n+j} \\
   & \s{3}{=} &
   \sum_{k}^{}\sum_{r}^{}
   (-1)^{n+p}\b{i-1}{k-1}\b{i-k}{r}\b{n-k-r}{p-k-r}\left(\sum_{j}^{}
\b{i-k-r}{j-k-r}(-a)^{j-k-r}\right)a^{r-k+2p-n} \\
   & \s{4}{=} &
   \sum_{k}^{}\sum_{r}(-1)^{n+p}\b{i-1}{k-1}\b{i-k}{r}\b{n-k-r}{p-k-r}
(-a)^{k+r-i}a^{r-k+2p-n} 
   \\
   & \s{5}{=} &
   \sum_{r}^{}(-1)^{n+p+r+i+1}\b{i-1}{r}\left(\sum_{k}^{}(-1)^{k-1}
   \b{i-r-1}{k-1}\b{n-k-r}{p-k-r}\right)a^{2r-i+2p-n} \\
   & \s{6}{=} & 
   \sum_{r=0}^{i-1}(-1)^{n+p+r+i+1}\b{i-1}{r}\b{n-i}{p-1-r}
       a^{2r-i+2p-n} \\
   & \s{7}{=} &
\sum_{k=1}^{i}(-1)^{n+p+k+i}\b{i-1}{k-1}\b{n-i}{p-k}a^{2k-2-i+2p-n} 
 \end{eqnarray*}
 and this last sum agrees with $\la_{p}u_{ip}$, as required.
 
 \noindent Notes:
 \vspace{-5mm}
 \begin{enumerate}
     \item  reverse sum over $j$ and reverse sum over $k$
 
     \item  apply $(*)$ with $N = n-k$ and $J = j-k$ to 
     $\b{n-j}{p-k}$, and 
     apply $(**)$ to $\b{i-1}{j-1}\b{j-1}{k-1}$
 
     \item  apply $(**)$ to $\b{i-k}{j-k}\b{j-k}{r}$ and rearrange sums
 
     \item  apply binomial theorem to sum over $j$, and use 
     $(1-a)^{i-k-r} = (-a)^{k+r-i}$
 
     \item  apply $(**)$ to $\b{i-1}{k-1}\b{i-k}{r}=\b{i-1}{i-k}\b{i-k}{r}$, and collect 
     terms in $k$
 
     \item  apply $(*)$ to evaluate sum over $k$

     \item  change summation index from $r$ to $k$ with $r = k-1$
 \end{enumerate}

The proof that $V^{2} = I_{n}$ uses the ``trinomial revision'' identity 
$(**)$ and the companion identity 
\[ 
\hspace{50mm}\b{I}{J}\b{J}{K} =
\b{I}{K}\b{I-K}{I-J}\hspace{20mm}\mathrm{integers\ }I,J,K\quad  (***)
\]
as well as the following three identities
\[
\hspace{40mm}\b{M}{k}\b{N}{L-k} = \b{M+N}{L}\hspace{35mm} 
\mathrm{integer\ }L\quad (\dag)
\]
\[
\hspace{50mm}\sum_{r}^{}(-1)^{r}\b{N}{r} =
\delta_{N,0}\hspace{30mm}\mathrm{integer\ }
N\ge 0\quad (\dag\dag)
\]
\[
\hspace{40mm}\sum_{u}^{}(-1)^{u}\b{N}{L-u}\b{N-L+u}{u} =
\delta_{L,0}\hspace{20mm} \mathrm{integer\ }L\quad (\dag\dag\dag)
\]
Identity (\dag) is the Vandermonde convolution, (\dag\dag) follows 
from $(1-1)^{N} = \delta_{N,0}$, (\dag\dag\dag) can be reduced to 
(\dag\dag) for $N \ge 0$ using (**), and holds for all $N$ since both 
sides are polynomials in $N$. 
Note that we don't need $a^{2} = a+1$: $V^{2} = I_{n}$ 
holds considering $V$ as a matrix with polynomial entries.

The identity $V^{2} = I_{n}$ is equivalent to 
$W^{2} = (1+a^{2})^{n-1}I_{n}$ with 
$w_{ij} = (-1)^{j}a^{n-j} \times $ 
\linebreak 
$\sum_{r}^{}(-1)^{i-r}\b{i-1}{r-1}\b{n-i}{j-r}a^{2r-i-1}$.
The equalities in the 
computation on the next page of the $(i,k)$ entry of $W^{2}$ are labelled 
with the identity used in that step or with a number referring to the 
following notes.

 \noindent Notes:
 \vspace{-5mm}
 \begin{enumerate}
     \item reverse sum on $j$

     \item apply the Vandermonde convolution to $\b{n-j}{s-1}$ and to 
           $\b{j-1}{k-s}$

     \item apply $(***)$ successively to rewrite the first 
           three factors, use the binomial theorem to evaluate the 
           parenthesized sum, and rearrange sums

     \item eliminate the sum on $s$

     \item rearrange sums

\end{enumerate}

\pagebreak

\begin{eqnarray*}	
 &  & \sum_{j=1}^{n}w_{ij}w_{jk}  \\
 & = & \sum_{j=1}^{n}\sum_{r,s}^{}(-1)^{r+s+i+k}\b{i-1}{r-1}
\b{j-1}{s-1}\b{n-i}{j-r}\b{n-j}{k-s}a^{2n+2r+2s-2j-i-k-2}  \\
& \s{1}{=} & \sum_{j=1}^{n}\sum_{r,s}^{}(-1)^{r+s+i+k}\b{i-1}{r-1}
\b{n-j}{s-1}\b{n-i}{n-j-r+1}\b{j-1}{k-s}a^{2r+2s+2j-i-k-4}  \\
 & \s{2}{=} &  \sum_{j,r,s,t,u}^{}(-1)^{r+s+i+k}\b{i-1}{r-1}
\b{r-1}{t-1}\b{n-j-r+1}{s-t}\b{n-i}{n-j-r+1}\b{j+r-i-1}{u} \times \\
& & \hspace*{40mm} \b{i-r}{k-s-u} a^{2r+2s+2j-i-k-4}  \\
	 &\s{(***)}{=} & \sum_{j,r,s,t,u}^{}(-1)^{r+s+i+k}\b{i-1}{r-1}
\b{r-1}{t-1}\b{n-i}{s-t}\b{n-i-s+t}{j+r-i-1}\b{j+r-i-1}{u} \times \\ 
& & \hspace*{40mm} \b{i-r}{k-s-u} a^{2r+2s+2j-i-k-4} \\
	 & \s{(**)}{=} & \sum_{r,s,t,u}^{}(-1)^{r+s+i+k}\b{i-1}{r-1}\b{r-1}{t-1} 
\b{i-r}{k-s-u}\b{n-i}{s-t}\b{n-i-s+t}{u} \times \\
& &\hspace*{40mm} \left(\sum_{j}^{}\b{n-i-s+t-u}{j+r-i-1-u}
        a^{2(j+r-i-u-1)}\right) a^{2s+i-k+2u-2} \\
	 & \s{3}{=} & \sum_{s,t,u}^{}(-1)^{s+i+k}\b{i-1}{t-1}\b{i-t}{k-s-u}
\left(\sum_{r}^{}(-1)^{r}\b{i-t-k+s+u}{r-t}\right) \times \\ 
& &\hspace*{40mm} \b{n-i}{s-t}\b{n-i-s+t}{u}(1+a^{2})^{n-i-s+t-u}
      a^{2s+i-k+2u-2}  \\
&\s{(\dag\dag)}{=} & \sum_{\substack{s,t,u\,:\,k-s-u \ge 0 \\ i-t=k-s-u}}(-1)^{s+i+k+t}
\b{i-1}{t-1}\b{i-t}{k-s-u} \times \\
& & \hspace*{40mm} \b{n-i}{s-t}\b{n-i-s+t}{u}(1+a^{2})^{n-i-s+t-u}
a^{2s+i-k+2u-2} \\
	 & \s{4}{=} &\sum_{t,u\,:\,t \le 
i}^{}(-1)^{u}\b{i-1}{t-1}\b{n-i}{k-i-u}\b{n+u-k}{u}(1+a^{2})^{n-k}
a^{k+2t-i-2}  \\
&\s{5}{=}  &(1+a^{2})^{n-k}a^{k-i}\sum_{t}^{}\b{i-1}{t-1}a^{2(t-1)}\sum_{u}(-1)^{u}
\b{n-i}{k-i-u}\b{n+u-k}{u}   \\
	 &\s{(\dag\dag\dag)}{=}  &(1+a^{2})^{n-k}a^{k-i}(1+a^{2})^{i-1}\delta_{i,k}   \\
	 & = & (1+a^{2})^{n-1}\delta_{i,k},\quad \mathrm{as\ required.} \\
\end{eqnarray*}

Finally, a mild generalization. Let $x$ be an indeterminate and let $R(x)$ 
be the $n \times n$ matrix with $(i,j)$ entry $\binom{i-1}{n-j}x^{i+j-n-1}$. Then the 
eigenvalues and eigenvectors of $R(x)$ are precisely as above but 
with $a$ a root of $a^{2} = ax+1$ rather than of $a^{2} = a+1$.

\end{document}